\documentclass[journal]{IEEEtran}
\ifCLASSINFOpdf
\else
\fi
\usepackage{amsmath}
\usepackage{amssymb}

\usepackage{graphicx}
\usepackage{subfigure}
\usepackage{amsfonts}
\usepackage{mathtools}
\usepackage{amsthm}
\usepackage[usenames]{color}
\usepackage[noabbrev]{cleveref}
\usepackage{psfrag}

\newtheorem{theorem}{Theorem}

\newtheorem{definition}[theorem]{Definition}
\newtheorem{example}[theorem]{Example}

\newtheorem{lemma}[theorem]{Lemma}

\newtheorem{remark}[theorem]{Remark}

\usepackage{algorithmic}

\usepackage{array}
\begin{document}
%
\title{Finite Necessary and Sufficient Stability Conditions for Linear System with Pointwise and Distributed Delays}
%
%
%

\author{Alejandro~Castaño,
        Carlos~Cuvas,
        Alexey~Egorov~\IEEEmembership{}%
        and~Sabine~Mondi\'e~\IEEEmembership{}
\thanks{A. Castaño and S. Mondi\'e are with the Department of Automatic Control, CINVESTAV-IPN, Mexico City, 07360, Mexico (e-mail: acastano@ctrl.cinvestav.mx; smondie@ctrl.cinvestav.mx).}
\thanks{C. Cuvas is with the Universidad Autónoma del Estado de Hidalgo, Pachuca de Soto 42000, Pachuca de Soto, Mexico (e-mail: carlos$\_$cuvas@uaeh.edu.mx).}
\thanks{A. Egorov is with St. Petersburg State
	University, St. Petersburg, 199034 Russia. e-mail: alexey.egorov@spbu.ru.}
}

%
%

\markboth{}%
{Castaño \MakeLowercase{\textit{et al.}}: Finite Necessary and Sufficient Stability Conditions for Linear System with Pointwise and Distributed Delays}
%



\maketitle

\begin{abstract}
This contribution presents two exponential stability criteria for linear systems with multiple pointwise and distributed delays. These results (necessary and sufficient conditions) are given in terms of the delay Lyapunov matrix and the fundamental matrix of the system. An important property is that the stability test requires a finite number of mathematical operations.
\end{abstract}


%
\IEEEpeerreviewmaketitle

\section{Introduction}
\label{sec:introduction}
Systems with distributed delay are those whose evolution depends not only on the present time dynamics but also on the dynamics at a previous time instant, as well as on the cumulative effect of past values of the dynamics \cite{kharitonov2013time}. These systems are used to model the time lag phenomenon in thermodynamics \cite{zheng2002robust}, population dynamics \cite{kolmanovskii2012applied}, traffic flow models \cite{juarez2020stability}, networked control systems \cite{moruarescu2007stability}, PID controller design \cite{kolmanovskii2012applied}, hematopoietic cell maturation \cite{ozbay2008stability}, among others. A distributed kernel allows thinner modeling of the interactions between the different system components. Consequently, due to the relevance, complexity, and potential applications, the stability study of systems with distributed delays has been of interest in recent decades, also representing a challenge from a purely theoretical perspective.

We consider linear time-delay systems of the form
\begin{equation}\label{Eq. System with distributed delay}
\dot{x}(t) = \sum_{j=0}^{m} A_{j} x(t-h_{j}) + \int_{-H}^{0} G(\theta) x(t+\theta) d\theta, \; \forall \; t \geq 0,
\end{equation}
where $A_{j}$, $j=0,\dots,m$, are real $n \times n$ matrices and $0 = h_{0} < h_{1} < \cdots < h_{m} = H$ are the delays. The function $G(\theta)$ is a real piecewise continuous matrix function defined for $\theta \in [-H,0]$ that represents the kernel of the distributed delay. 

The time-domain stability analysis is mainly based on the ideas introduced in \cite{krasovskii1956application}, extending Lyapunov's classical approach for delay-free systems to the time-delay case. This approach replaces classical Lyapunov functions that depend on the instantaneous state of a system with functionals that depend on the state defined on the delay interval. The converse approach for the construction of this class of functionals was originally presented in \cite{repin1965quadratic}, \cite{datko1972algorithm} and \cite{wenzhang1989generalization}. In \cite{kharitonov2003lyapunov}, Lyapunov-Krasovskii functionals of complete type, admitting quadratic lower bounds when the system is exponentially stable, were proposed for pointwise delay systems. These functionals are determined by a matrix function called delay Lyapunov matrix defined on the delay interval, which satisfies four properties:  continuity, dynamic, symmetric, and algebraic \cite{kharitonov2013time}. 

The study of complete type functionals is not restricted to pointwise delay systems, as it also covers systems with distributed delays \cite{kharitonov2013time}. Because of the distributed nature of the delay, the Lyapunov-Krasovskii functional with prescribed negative quadratic derivative now involves double and triple integral terms \cite{kharitonov2013time}. The case of general kernels has been investigated by \cite{solomon2013new} where the application of new integral inequalities is suggested. Complete type functionals have been used in the context of robust stability \cite{gouaisbaut2015stability}, exponential stability analysis, \cite{egorov2017necessary}, and in determining the parameters or delays critical values  \cite{ochoa2009time}, just to mention a few. In these contribution the system was assumed to be stable.

In recent years, the converse approach allowed revisiting the plain time-delay stability problem from the perspective of necessity. The goal is to extend the exponential stability criterion of delay-free systems stated in terms of the positivity of the Lyapunov matrix $P$, solution of the Lyapunov equation $A^{T}P + PA = -Q$. This was first successfully addressed in \cite{mondie2012assessing} and \cite{egorov2013stability} for the case of scalar delay equations. Families of necessary stability conditions were obtained for linear systems with multiple delays \cite{egorov2014necessaryb}, \cite{egorov2014bnecessary}, neutral-type delays \cite{gomez2019necessary}, and distributed delays \cite{cuvas2015necessary}. \textcolor{black}{ Many examples in these classes indicated the condition might be sufficient. Sufficiency in an infinite number of operations was established in \cite{egorov2014new}, and \cite{egorov2017necessary}. Finiteness of the criteria, a highly valuable feature both from theoretical and practical points of view}  was achieved for systems with multiple pointwise delays \cite{gomez2019lyapunov}, and of neutral type \cite{gomez2020necessary}, but for systems with distributed delays, it is an open problem.

%
%

%
%
The main contribution of this manuscript is to present two finite exponential stability conditions for systems of the form \eqref{Eq. System with distributed delay}. The first one is exclusively given in terms of the delay Lyapunov matrix and uses arguments similar to those used in \cite{gomez2017lyapunov} and \cite{gomez2020necessary}. The second depends also on the system's fundamental matrix but consists of a reduced number of mathematical operations.

The paper’s organization is as follows: Section \ref{Sec: Preliminaries1} is dedicated to some preliminaries on the systems with distributed delays. In Section \ref{sec: Lyapunov-Krasovskii framework}, results on Lyapunov-Krasovskii functionals with prescribed derivatives are recalled, and some stability theorems for system \eqref{Eq. System with distributed delay} are introduced. The choice of initial functions leading to an expression of the functional in terms of the delay Lyapunov matrix, without the assumption of the stability of the system, is discussed in Section \ref{sec:instrumental result}. The main results are given in Section \ref{sec:Main result}. The paper ends with illustrative examples in Section \ref{sec: examples}, followed by some concluding remarks.

\emph{Notation:} We denote the space of piecewise continuous and continuously differentiable functions by $\mathcal{H} = PC([-H,0],\mathbb{R}^{n})$ and $\mathcal{C}^{(1)}([-H,0],\mathbb{R}^{n})$, respectively. For vectors and matrices we use the Euclidean norm, denoted by $\|\cdot\|$, and for functions $\varphi$, we use the uniform norm $\displaystyle \|\varphi\|_{H} = \sup_{\theta \in  [-H,0]}\|\varphi(\theta)\|.$
The transpose of a matrix $A$ is denoted by $A^T$, while the minimum and maximum eigenvalue of a symmetric matrix $Q$ are represented by $\lambda_{\min}(Q)$ and $\lambda_{\max}(Q)$,  respectively. The notation $Q>0$ means that the symmetric matrix $Q$ is positive definite. The symbol $\lceil \cdot \rceil$ denotes the ceiling function. The identity $q \times q$ matrix is denoted by $I_{q}$.

\section{Preliminaries} 
\label{Sec: Preliminaries1}
In this section, essential concepts in the analysis in the time-domain \cite{kharitonov2013time} of linear systems with pointwise and distributed delays are recalled.

\subsection{System}
Some basic definitions on system \eqref{Eq. System with distributed delay} are first introduced. The initial function $\varphi$ is \textcolor{black}{assumed} to be piecewise continuous, $\varphi \in \mathcal{H}$. The restriction of the solution $x(t,\varphi)$ of system \eqref{Eq. System with distributed delay} on the segment $[t-H,t]$ is defined by
\begin{equation*}
x_{t}(\varphi) : \theta \; \rightarrow \; x(t+\theta,\varphi), \quad \theta \in [-H,0].
\end{equation*}

\begin{definition}[see \cite{bellman1963differential}]
\label{def:ExponentialStability}
System \eqref{Eq. System with distributed delay} is said to be exponentially stable if there exist $\gamma \geq 1$ and $\sigma > 0$ such that 
\begin{equation*}
    \| x(t,\varphi) \| \leq \gamma e^{\sigma t} \| \varphi \|_{H}, \quad t \geq 0.
\end{equation*}
\end{definition}

The fundamental matrix of system \eqref{Eq. System with distributed delay}, denoted by $K(t)$, satisfies \cite{bellman1963differential}:
\begin{equation*}
\dot{K}(t) = \sum_{j=0}^{m}  K(t-h_{j}) A_{j} + \int_{-H}^{0}  K(t+\theta) G(\theta)  d\theta, \quad t \geq 0,
\end{equation*} 
with the initial condition $K(t) = 0$ for $t < 0$ and $K(0) = I_{n}$.
\begin{remark}
\label{rem:fundamental matrix}
The fundamental matrix also satisfies the matrix equation 
\begin{equation}\label{Eq. Fundamental matrix distributed delay}
\dot{K}(t) = \sum_{j=0}^{m} A_{j} K(t-h_{j}) + \int_{-H}^{0} G(\theta) K(t+\theta) d\theta.
\end{equation} 
\end{remark}
\textcolor{black}{However, the products of the matrices $A_{j}$, $j=0,1,\dots,m$ with the fundamental matrix, do not commute individually.} 
\begin{lemma}
\label{rem: fundamental matrix bound}
The matrix $K(t)$ satisfies
\begin{equation*}
\| K(t) \| \leq e^{M_{1} t} 
\end{equation*}
where $M_{1} = M+ b H$, with $b = \underset{\theta \in[-H,0]}{\sup}\| G(\theta)\|$ and $M = \sum_{j=0}^{m}\| A_{j}\| $. Hence, one can calculate the number $L$, such that
\begin{equation*}
\begin{split}
\left\| \dfrac{d {K}(t)}{dt} \right\| &\leq  M_{2} e^{M_{1}t} =L, \quad t\in[0,H]
\end{split}
\end{equation*}
where $\displaystyle M_{2} =  \sum_{j=0}^{m} \|A_{j}\|  e^{-M_{1}h_{j}} + \frac{ b  }{M_{1} }  (1-e^{-M_{1} H})$.
\end{lemma} 
\begin{proof}
\textcolor{black}{Integrating the expression \eqref{Eq. Fundamental matrix distributed delay} from $0$ to $t$, and applying the Cauchy-Schwarz inequality, we have}
\begin{equation*}
\begin{split}
\| K(t) \| &\leq  \| K(0) \| +  (M + b H) \int_{0}^{t} \|K(\tau)\|d\tau \\
          &\leq  1 +  M_{1} \int_{0}^{t} \|K(\tau)\|d\tau. \\
\end{split}
\end{equation*}
\textcolor{black}{Now, applying the Bellman-Gronwall lemma to the previous expression, we obtain}
\begin{equation*}
\begin{split}
\| K(t) \| &\leq  e^{M_{1} t}. \\
\end{split}
\end{equation*}
\textcolor{black}{Finally, the last inequality and equation \eqref{Eq. Fundamental matrix distributed delay}imply that}
\begin{equation*}
\begin{split}
\| \dot{K}(t) \| &\leq \sum_{j=0}^{m} \|A_{j}\|  e^{M_{1} (t-h_{j})}  + b \int_{-H}^{0}  e^{M_{1} (t+\theta)}  d\theta\\
                 &\leq  M_{2} e^{M_{1}t}, 
\end{split}
\end{equation*}
\textcolor{black}{with $M_{2}$ expressed in Lemma \ref{rem: fundamental matrix bound}, which completes the proof.}
\end{proof}


\section{Lyapunov-Krasovskii framework}
\label{sec: Lyapunov-Krasovskii framework}
The functional $v_{0}(x_{t}(\varphi))$ with prescribed derivative along the trajectories of system \eqref{Eq. System with distributed delay}, given by
\begin{equation*}
\dfrac{d v_{0}(x_{t}(\varphi))}{dt} = -x^{T}(t,\varphi)Wx(t,\varphi),
\end{equation*} 
where $W$ is a positive definite matrix, was introduced in \cite{wenzhang1989generalization}. It has the form 
\begin{equation*}\label{Eq. Quadratic functional v0}
\begin{split}
v_{0}(\varphi) &= \varphi^{T}(0) U(0) \varphi(0) 
               + 2 \varphi^{T}(0) \sum_{j=1}^{m}\int_{-h_{j}}^{0}U(-\theta-h_{j})\\ 
               &\times A_{j}\varphi(\theta)d\theta
               + \sum_{k=1}^{m} \sum_{j=1}^{m} \int_{-h_{k}}^{0}\varphi^{T}(\theta_{1})A_{k}^{T} \\
                & \times \int_{-h_{j}}^{0} U(\theta_{1}+h_{k}-\theta_{2}-h_{j})A_{j}\varphi(\theta_{2})d\theta_{2}d\theta_{1}\\
                 &+ 2 \varphi^{T}(0) \int_{-H}^{0} \int_{-H}^{\theta}U(\xi-\theta)G(\xi)d\xi\varphi(\theta) d\theta\\
               \end{split}
\end{equation*}
\begin{equation*}
\begin{split}
                &+ 2 \sum_{j=1}^{m} \int_{-h_{j}}^{0} \int_{-H}^{0} \int_{-H}^{\theta_{2}} \varphi^{T}(\theta_{1})A_{j}^{T}\\
                      & \times U(h_{j}+\theta_{1}-\theta_{2}+\xi)G(\xi)\varphi(\theta_{2})d\xi d\theta_{2} d\theta_{1}\\
                      &+ \int_{-H}^{0} \varphi^{T}(\theta_{1}) \int_{-H}^{0} \int_{-H}^{\theta_{1}} \int_{-H}^{\theta_{2}}G^{T}(\xi_{1}) \\ 
                      & \times  U(\theta_{1}-\theta_{2}-\xi_{1}+\xi_{2})G(\xi_{2})d\xi_{2}d\xi_{1} \varphi(\theta_{2})d\theta_{2} d\theta_{1}.
\end{split}
\end{equation*}
Each term of the foregoing equation depends on the matrix-valued function $U(\tau)$, named the \textit{delay Lyapunov matrix} associated with $W$. It satisfies the following set of properties,
\begin{enumerate}
    \item Continuity property
    \begin{equation*}
        U\in \mathcal{C}\left(\mathbb{R},\mathbb{R}^{n\times n}\right).
    \end{equation*} 
    \item Dynamic property
    \begin{equation*}\label{Eq. dynamic property}
        \dfrac{d}{d\tau} U(\tau) = \sum_{j=0}^{m}U(\tau-h_{j})A_{j}+\int_{-H}^{0}U(\tau+\theta)G(\theta)d\theta. 
    \end{equation*}
    \item Symmetry property
    \begin{equation*}\label{Eq. symmetry property}
     U(\tau) = U^{T}(-\tau) , \quad \tau \geq 0. 
    \end{equation*}
    \item Algebraic property
    \begin{equation*}\label{Eq. algebraic property}
     \begin{split}
         -W &= \sum_{j=0}^{m}[A_{j}^{T}U^{T}(-h_{j})+U(-h_{j})A_{j}] \\ &+ \int_{-H}^{0}[G^{T}(\theta)U^{T}(\theta)+U(\theta)G(\theta)]d\theta. 
     \end{split}   
    \end{equation*}
\end{enumerate}
The uniqueness of the delay Lyapunov matrix is established in the next theorem. 
\begin{theorem}[see \cite{kharitonov2013time}]\label{Theorem Lyapunov Condition}
System \eqref{Eq. System with distributed delay} admits a unique delay Lyapunov matrix $U$ associated with a matrix $W$ if and only if the system satisfies the Lyapunov condition, i.e., if the spectrum 
\begin{equation*}
\Lambda = \left\{s \; | \; \textrm{det}\left( sI-\sum_{j=0}^{m}A_{j}e^{-sh_{j}} - \int_{-H}^{0}e^{s \theta}G(\theta)d\theta \right)=0  \right\}
\end{equation*}
does not contain any root $\hat{s}$ such that $-\hat{s}$ is also a root.
\end{theorem}
Let us introduce now the following quadratic functional $v_{1}$
\begin{equation*}\label{Eq. Functional v1(phi)}
v_{1}(\varphi) = v_{0}(\varphi) + \int_{-H}^{0} \varphi^{T}(\theta) W \varphi(\theta) d\theta,
\end{equation*}
whose derivative along the solution of system \eqref{Eq. System with distributed delay} is
\begin{equation}\label{Eq. Derivative of v1(t)}
\dfrac{d v_{1}(x_{t}(\varphi))}{dt} = -x^{T}(t-H,\varphi) W x(t-H,\varphi).
\end{equation}
This functional is not of complete type, i.e., its derivative does not include the whole delayed state, a valuable property for robust stability analysis. However, it allows presenting two significant stability/instability results. The first one is the existence of a quadratic lower bound of $v_{1}$ when the system is stable.
\begin{theorem}\label{Theorem Stability criterion v1>a1}
If system \eqref{Eq. System with distributed delay} is exponentially stable, then there exist positive numbers $\alpha_{0}$ and $\alpha_{1}$ such that for any $\varphi \in \mathcal{H}$
\begin{equation} \label{Eq. v1 > a |phi|}
v_{1}(\varphi) \geq \alpha_{0}\|\varphi(0)\|^{2} + \alpha_{1} \int_{-H}^{0} \|\varphi(\theta) \|^{2} d\theta.
\end{equation}
\end{theorem}
\begin{proof}
We define an auxiliary functional of the form
\begin{equation*}
\begin{split}
        \Tilde{v}_{1}(\varphi) &= v_{1}(\varphi) - \dfrac{1}{m+1} \sum_{j=1}^{m} \int_{-H}^{-h_{j}}\varphi^{T}(\theta)W \varphi(\theta)d\theta \\ & - \dfrac{1}{m+1} \int_{-H}^{0}\varphi^{T}(\theta)\left[ W - \dfrac{H+\theta}{\eta H}W \right] \varphi(\theta)d\theta \\ &- \alpha_{0} \|x(t)\|^{2},
\end{split}
\end{equation*}
where $\alpha_{0}$ is assumed to be a positive constant, and $\eta > 1$.  The time derivative of this auxiliary functional along the solution of system \eqref{Eq. System with distributed delay} is
\begin{equation*}
    \begin{split}
        \left.\dfrac{d}{dt}\Tilde{v}_{1}(x_{t}) \right|_{\eqref{Eq. System with distributed delay}} &= -\Tilde{w}(x_{t}),
    \end{split}
\end{equation*}
where
\begin{equation*}
    \begin{split}
        \Tilde{w}(x_{t}) =&  
         \dfrac{1}{m+1} \sum_{j=0}^{m} x^{T}(t-h_{j}) W x(t-h_{j})\\ 
         -& \dfrac{1}{\eta(m+1)} x^{T}(t) W x(t) 
         +2\alpha_{0} x^{T}(t) \sum_{j=0}^{m} A_{j} x(t-h_{j}) \\ 
         +&  \dfrac{1}{\eta H(m+1)} \int_{-H}^{0} x^{T}(t+\theta) W x(t+\theta)d\theta\\
         +& 2\alpha_{0} x^{T}(t) \int_{-H}^{0} G(\theta) x(t+\theta) d\theta, \quad t \geq 0.\\
    \end{split}
\end{equation*}
The inequality $2a^{T}b \geq -a^{T}a - b^{T}b$ implies that the functional $\Tilde{w}(\varphi)$ admits a lower estimation of the form
\begin{equation*}
    \Tilde{w}(\varphi) \geq \Tilde{\varphi}^{T}R_{1}(\alpha_{0})\Tilde{\varphi} + \int_{-H}^{0}\varphi^{T}(\theta)R_{2}(\theta,\alpha_{0}) \varphi(\theta)d\theta,
\end{equation*}
where
\begin{equation*}
    \begin{split}
    \Tilde{\varphi} &= (\varphi^{T}(0) , \varphi^{T}(-h_{1}) , \dots , \varphi^{T}(-H))^{T},\\
        R_{1}(\alpha_{0})&=\dfrac{1}{m+1}\left[ I_{m+1} \otimes W -  \dfrac{1}{\eta} (e_{1}\cdot e^{T}_{1}) \otimes W \right]+ \alpha_{0} \\ &\cdot \left[ e_{1}\otimes I_{n} \cdot A + A^{T} \cdot e_{1} \otimes I_{n} -H \cdot (e_{1}\cdot e^{T}_{1}) \otimes I_{n} \right],\\
        R_{2}(\alpha_{0}) &= \dfrac{1}{\eta H(m+1)} W - \alpha_{0} G^{T}(\theta) G(\theta),
    \end{split}
\end{equation*}
and the block-matrix $A=(A_{0},A_{1},\dots,A_{m}) \in \mathbb{R}^{n \times n (m+1)}$ and $e_{1} = (1,0,\dots,0)^{T} \in \mathbb{R}^{m+1}$. The matrix $W$ is positive definite, so there exists $\alpha_{0}>0$ such that the following inequalities hold
\begin{equation*}
    \begin{split}
        R_{1}(\alpha_{0}) \geq 0, \quad R_{2}(\theta,\alpha_{0}) \geq 0, \quad \theta \in [-H,0].
    \end{split}
\end{equation*}
Therefore, by choosing $\alpha_{0} =  \alpha_{0}^{\star}$ with
\begin{equation}\label{Eq. value of alpha0*}
        \alpha_{0}^{\star} = \min(\alpha_{01}^{\star},\alpha_{02}^{\star}),\\
\end{equation}
where
\begin{equation*}
    \begin{split}
        \alpha_{01}^{\star} &= - \dfrac{1}{(m+1)\lambda_{min}(P_{1})},\quad 
        \alpha_{02}^{\star} = \dfrac{1}{\eta H (m+1)\lambda_{max}(P_{2})},\\
        P_{1} &= \left[ I_{m+1} \otimes W -  \dfrac{1}{\eta} (e_{1}\cdot e^{T}_{1}) \otimes W \right]^{-1}\\ &\cdot \left[ e_{1}\otimes I_{n} \cdot A + A^{T} \cdot e_{1} \otimes I_{n} -H \cdot (e_{1}\cdot e^{T}_{1}) \otimes I_{n} \right]\\
        P_{2} &=  W^{-1}G^{T}(\theta)G(\theta), \quad \forall \; \theta \in [-H,0],\\
    \end{split}
\end{equation*}
we have that $R_{1}(\alpha_{0}) \geq 0$ and $R_{2}(\theta,\alpha_{0}) \geq 0$. For the such value of $\alpha_{0}$ we obtain
\begin{equation*}
    \Tilde{w}(\varphi) \geq 0, \quad \varphi \in PC([-H,0],\mathbb{R}^{n}).
\end{equation*}
As the system is exponentially stable, i.e., $\displaystyle  \lim_{t \to \infty} \Tilde{v}_{1}(x_{t})=0$, then
\begin{equation*}
\begin{split}
 \Tilde{v}_{1}(\varphi)
 &=\int_{0}^{\infty}\Tilde{w}(x_{t}(\varphi)) dt \geq 0.
 \end{split}
\end{equation*}
Therefore
\begin{equation*}
\begin{split}
       v_{1}(\varphi) &\geq \dfrac{1}{m+1} \sum_{j=1}^{m} \int_{-H}^{-h_{j}}\varphi^{T}(\theta)W \varphi(\theta)d\theta + \alpha_{0} \|x(t)\|^{2} \\ & + \dfrac{1}{m+1} \int_{-H}^{0}\varphi^{T}(\theta)\left[ W - \dfrac{H+\theta}{\eta H}W \right] \varphi(\theta)d\theta \\
       & \geq \alpha_{0} \|x(t)\|^{2}  + \dfrac{1}{m+1} \int_{-H}^{0}\varphi^{T}(\theta)W_{1} \varphi(\theta)d\theta. 
\end{split}
\end{equation*}
The result follows by setting $\alpha_{1} = \frac{\lambda_{min}(W_{1})}{m+1}$, with $W_{1}=\left[ W - \frac{H+\theta}{\eta H}W \right]$.
\end{proof}
The second result concerns instability. It is based on the ideas introduced in \cite{medvedeva2013constructive}, that establish that if the system is unstable, the functional $v_{1}$ does not admit a positive lower bound.
\begin{theorem}[see \cite{egorov2017necessary}]\label{Theorem v1 <= -alpha1}
If system \eqref{Eq. System with distributed delay} is unstable and satisfies the Lyapunov condition, then for every $\hat{\alpha}_{1} > 0$ there exists a function $\hat{\varphi} \in \mathcal{H}$ such that
\begin{equation*}
v_{1}(\hat{\varphi}) \leq - \hat{\alpha}_{1}.
\end{equation*}  
\end{theorem}
We recall the following bilinear functional \cite{cuvas2015necessary},
\begin{equation}\label{Eq. Bilinear functional z}
z(\varphi_{1},\varphi_{2}) = \dfrac{1}{4} \left( v_{1}(\varphi_{1} + \varphi_{2}) - v_{1}(\varphi_{1} - \varphi_{2}) \right)\\
\end{equation}
where $\varphi_{1},\varphi_{2} \in \mathcal{H}$. This functional reduces to the functional $v_{1}(\varphi)$ when $\varphi_{1} = \varphi_{2}=\varphi$, i.e., $v_{1}(\varphi)=z(\varphi,\varphi)$. In the next lemma, we give an upper bound for the functionals $v_{1}$ and $z$. 
\begin{lemma}\label{Lemma |v1|<alpha2 |z|<alpha2}
For any $\varphi$, $\varphi_{1}$, $\varphi_{2}$ $\in \mathcal{H}$, there exists $\alpha_{2}>0$ such that
\begin{equation*}
\begin{split}
|z(\varphi_{1},\varphi_{2})| &\leq \alpha_{2} \|\varphi_{1}\|_{H} \|\varphi_{2}\|_{H},\\
|v_{1}(\varphi)| &\leq \alpha_{2} \|\varphi\|_{H}^{2},\\
\end{split}
\end{equation*}
where
\begin{equation*}
\begin{split}
\alpha_{2} &= \nu \left( 1 +  N_{1} \right)^{2}  + b \nu H^{2}\left( 1 + N_{1} \right) + \dfrac{1}{4} b^{2} \nu H^{2} + H \| W \|,
\end{split}
\end{equation*}
with
\begin{equation*}
\begin{split}
\nu = \max_{\tau\in [0,H]}\|U(\tau)\|, \;
b = \sup_{\theta \in[-H,0]}\| G(\theta)\|, \; N_{1} = \sum_{j=1}^{m} \| A_{j} \| h_{j}.
\end{split}
\end{equation*}
\end{lemma}
\begin{proof}
The first inequality can be deduced by applying the Cauchy-Schwarz inequality to each term of the functional $z(\varphi_{1},\varphi_{2})$ defined in \eqref{Eq. Bilinear functional z}. The second inequality follows from the fact that $v_{1}(\varphi)=z(\varphi,\varphi)$.
\end{proof}


\section{Instrumental results}
\label{sec:instrumental result}
We introduce some key auxiliary results that will be crucial for the proof of the main theorems in Section \ref{sec:Main result}. An important element in the proof of the stability criterion presented by \cite{egorov2016finite} is the following compact set in the space of continuously differentiable functions:
\begin{equation}\label{Eq. Compact set S}
\begin{split}
\mathcal{S} &= \{ \varphi \in \mathcal{C}^{(1)}([-H,0],\mathbb{R}^{n}) \; :\\ & \; \| \varphi \|_{H} = \|\varphi(0)\|=1 \; ; \; \| \varphi(\theta) \| \leq 1, \; \|\varphi'\|\leq M_{1} \},
\end{split}
\end{equation}
with $M_{1} = M + b$. Consider the function $\psi_{r}:[-H,0] \to \mathbb{R}^{n}$ defined as
\begin{equation}\label{Eq. Function psi_r}
\psi_{r}(\theta) = \sum_{i=1}^{r}K(\tau_{i}+\theta)\gamma_{i}, \quad \theta \in [-H,0],
\end{equation}
where $r \in \mathbb{N}$, $\tau_{i} \in [0,H]$, and $\gamma_{i} \in \mathbb{R}^{n}$, $i=\overline{1,r}$ are arbitrary vectors. In \cite{cuvas2015necessary} and \cite{egorov2017necessary}, by introducing new properties that \textcolor{black}{connect} the delay Lyapunov matrix $U$ with the fundamental matrix, and using the function $\psi_{r}$ and bilinear functional \eqref{Eq. Bilinear functional z}, it is shown that
\begin{equation}\label{Eq. v1(psi_r)}
v_{1}(\psi_{r}) = \gamma^{T} \{ U(-\tau_{i}+\tau_{j}) \}_{i,j=1}^{r} \gamma = \gamma^{T} \mathcal{K}_{r}(\tau_{1},\dots,\tau_{r}) \gamma ,
\end{equation} 
where $\gamma = (\gamma_{1}^{T},\dots,\gamma_{r}^{T})^{T}$ and $\mathcal{K}_{r} \in \mathbb{R}^{n r \times n r}$. 
In \cite{egorov2017necessary}, it is shown that by appropriately choosing $\tau_{i}$ and $\gamma_{i}$, $i=\overline{1,r}$, any arbitrary function from the set $\mathcal{S}$ can be approximated by a function of the form \eqref{Eq. Function psi_r}. In the case of equidistant points in the interval $[0, H]$, namely,
\begin{equation*}
\tau_{i} = \dfrac{i-1}{r-1}H,
\end{equation*} 
the matrix $\mathcal{K}_{r}$ takes the form
\begin{equation*}
\begin{split}
  \mathcal{K}_{r} &= \left( 0, \dfrac{H}{r-1}, \dots,\dfrac{(r-2)H}{r-1},H  \right)  
  = \left\{  U \left( \dfrac{j-i}{r-1} H \right)  \right\}_{i,j=1}^{r}
\end{split}
\end{equation*}
where $r \geq 2$, and for $r=1$, $\mathcal{K}_{1} := U(0)$. Based on the fact that every continuous function can be approximated by the function $\psi_{r}$ in \eqref{Eq. Function psi_r}, the following stability criterion was presented in \cite{egorov2017necessary}:
\begin{theorem}\label{Theorem Criterio Kr} System \eqref{Eq. System with distributed delay} is exponentially stable if and only if the Lyapunov condition holds and for every natural number $r\geq 2$ 
\begin{equation}
\mathcal{K}_{r}(\tau_{1},\dots,\tau_{r}) = \left\{ U \left( \dfrac{j-i}{r-1}H \right) \right\}_{i,j=1}^{r} > 0.
\end{equation}
Moreover, if the Lyapunov condition holds and system \eqref{Eq. System with distributed delay} is unstable, then there exists a natural number $r$ such that $\mathcal{K}_{r}(\tau_{1},\dots,\tau_{r}) \ngeq 0$.
\end{theorem}
The following crucial result allows us to obtain an estimate of the approximation error, $E_{r} = \varphi - \psi_{r}$, by considering that the vectors $\gamma_{i}$, $i=\overline{1,r}$, are such that 
\begin{equation*}
\psi_{r}(-\tau_{i}) = \varphi(-\tau_{i}), \quad i = \overline{1,r}.
\end{equation*}
\begin{lemma}[see \cite{egorov2017necessary}]\label{Lemma |Er|<epsilon_r} 
For any $\varphi \in \mathcal{S}$, there exists a function $\psi_{r}$ of the form \eqref{Eq. Function psi_r} such that
\begin{equation*}
\|\varphi - \psi_{r}\|_{H} \leq \varepsilon_{r}=\dfrac{H(L+M_{1})e^{LH}}{r-1 +LH}.
\end{equation*}
Here $L$ is such that $\|\dot{K}(t)\| \leq L$, \textcolor{black}{almost everywhere} on $[0,H]$. 
\end{lemma}
The following result strengthens Theorem \ref{Theorem Stability criterion v1>a1}.
\begin{theorem}\label{Theorem: v1>=alpha0estrella}
If system \eqref{Eq. System with distributed delay} is exponentially stable, then for any $\varphi \in \mathcal{S}$
\begin{equation}\label{Eq. lower bound v1>=alpha0*}
    v_{1}(\varphi) \geq \alpha_{0}^{\star},
\end{equation}
where $\alpha_{0}^{\star}$ is expressed by \eqref{Eq. value of alpha0*}.
\end{theorem}
\begin{proof}
The result can be achieved by following the steps of the proof of Theorem \ref{Theorem Stability criterion v1>a1}, with $\varphi \in \mathcal{S}$. 
\end{proof}
We next introduce an instability result. The basic idea is inspired by research works presented in \cite{medvedeva2013constructive}, \cite{cuvas2015necessary} and \cite{egorov2017necessary}. We first introduce one auxiliary result.
\begin{lemma}[see \cite{gomez2017lyapunov}] \label{Lema de det(F)=0} 
Let $F \in \mathbb{C}^{n \times n}$ be a complex matrix. If $\det(F)=0$, then there exist $C_{1}, C_{2} \in \mathbb{R}^{n} $ such that
\begin{enumerate} 
    \item $ F( C_{1} + j C_{2} ) = 0$,
    \item $\Vert C_{1}\Vert =1$,
    \item $\Vert C_{2}\Vert \leq 1$,
    \item $C_{1}^{T}C_{2}=0$.
\end{enumerate}
\end{lemma}
The following theorem provides a necessary instability condition for system \eqref{Eq. System with distributed delay} based on a computable upper bound of the functional $v_{1}$. This bound is the cornerstone of the main result presented in this contribution.
\begin{theorem}\label{Theorem Stability criterion v1<-a1 }
If system \eqref{Eq. System with distributed delay} has an eigenvalue with a strictly positive real part, there exists $\varphi \in \mathcal{S}$ such that 
\begin{equation}\label{Eq. v1<alpha1 varphi en S}
v_{1}(\varphi) \leq - \alpha_{1} = - \dfrac{\lambda_{min}(W)}{4M_{1}}e^{-2M_{1} H} \cos^{2}(b_{1}),
\end{equation}
where $b_{1}$ exists and is a unique solution of the equation 
\begin{equation}\label{Eq. Expression g(b)}
\sin ^{4}(b_{1})((HM_{1})^{2}+b_{1}^{2})-(HM_{1})^{2}=0, \quad b_{1} \in \left( 0, \dfrac{\pi}{2} \right).
\end{equation}
\end{theorem}
\begin{proof}
As system \eqref{Eq. System with distributed delay} is unstable, there exists an eigenvalue $s_{0} = \alpha + i\beta$ with $\alpha>0$ and $\beta \geq 0$, and two vectors $C_{1},C_{2} \in \mathbb{R}^{n}$ that satisfy the condition of Lemma \ref{Lema de det(F)=0} with 
\begin{equation*}
F = s_{0}I-\sum_{j=0}^{m}A_{j}e^{-s_{0}h_{j}} - \int_{-H}^{0}e^{s_{0} \theta}G(\theta)d\theta.
\end{equation*}
In this case, the following expression is a particular solution of system \eqref{Eq. System with distributed delay} on $t\in(-\infty,\infty)$,
\begin{equation*}
\bar{x}(t) = e^{\alpha t}\phi (t), \quad \phi (t) = \cos (\beta t)C_{1}-\sin (\beta t)C_{2},
\end{equation*} 
with initial condition
\begin{equation*}
\varphi(\theta) = \bar{x}(\theta,\varphi), \quad \theta \in [-H,0].
\end{equation*}
The first step is to prove that $\varphi \in \mathcal{S}$. By Lemma \ref{Lema de det(F)=0}, $\| \varphi(0) \|=1$ and $\|\phi(t)\|^{2}=\cos^{2}(\beta t)\|C_{1}\|^{2}+\sin^{2}(\beta t)\|C_{2}\|^{2} \leq 1$.
The last expression implies that $\underset{t\in\mathbb{R}}{\textrm{max}}\|\phi(t)\|=1$, hence,
\begin{equation*}
\| \bar{x}(t,\varphi)\| = e^{\alpha t}\|\phi(t)\| \leq \|\varphi(0)\|=1, \quad t\leq 0.
\end{equation*}
Now, since $\bar{x}(t,\varphi)$ satisfies \eqref{Eq. System with distributed delay} for $t \in (-\infty,\infty)$, we have
\begin{equation*}
\begin{split}
\|\dot{\varphi}(t)\| &= \| \dot{\bar{x}}(t,\varphi) \|  \leq\sum_{j=0}^{m}\| A_{j}\| \|\bar{x}(t-h_{j},\varphi)\| \\ &+ \int_{-H}^{0} \|G(\theta)\| \|\bar{x}(t+\theta)\| d\theta \leq M_{1}, \quad \forall \; t\in[-H,0].
\end{split}
\end{equation*}
We need to estimate  $\alpha$. The following equation holds:
\begin{equation*}
s_{0}I  = \sum_{j=0}^{m}A_{j}e^{-s_{0}h_{j}} + \int_{-H}^{0}e^{s_{0} \theta}G(\theta)d\theta,
\end{equation*}
therefore
\begin{equation*}
|s_{0}| \leq \sum_{j=0}^{m} \|A_{j}\| |e^{-s_{0}h_{j}}| + \int_{-H}^{0}|e^{s_{0} \theta}| \|G(\theta)\|d\theta,\\
\end{equation*}
hence
\begin{equation*}\label{Eq. Cota de alpha<M1}
\alpha  \leq |s_{0}| \leq  \sum_{j=0}^{m} \|A_{j}\|  + b H = M_{1}.
\end{equation*}

Integrating the derivative of $v_{1}$ defined by \eqref{Eq. Derivative of v1(t)} from $0$ to $T$, where $T=2\pi/\beta$, if $\beta \neq 0$, and $T=1$, if $\beta = 0$, we obtain
\begin{equation*}
v_{1}(\bar{x}_{0})=v_{1}(\bar{x}_{T})  + \int_{-H}^{T-H}\bar{x}^{T}(t)W\bar{x}%
(t)dt,
\end{equation*}
where $W$ is positive definite. Since $T$ is the period of the function $\phi(t)$ for $t \in (-\infty,\infty)$, we have $\bar{x}(t+T)=e^{\alpha T}\bar{x}(t)$, and $v_{1}(\bar{x}_{T})=e^{2\alpha t}v_{1}(\bar{x}_{0})$, which implies that 
\begin{equation*}
\begin{split}
v_{1}(\bar{x}_{0}) &=-\dfrac{1}{e^{2\alpha T}-1}\int_{-H}^{T-H}\bar{x}^{T}(t)W\bar{%
x}(t)dt\\ &\leq -\dfrac{\lambda _{min}(W)}{e^{2\alpha T}-1}\int_{-H}^{T-H}\Vert 
\bar{x}(t)\Vert ^{2}dt.
\label{intermed1}
\end{split}
\end{equation*}
Substituting the particular solution $\bar{x}(t)$ in the integral term of the previous expression, we get
\begin{equation*}
\begin{split}
&\int_{-H}^{T-H}\Vert \bar{x}(t)\Vert ^{2}dt = \int_{-H}^{T-H}e^{2\alpha t} \| \phi(t) \|^{2} dt\\
&\geq \int_{-H}^{T-H}e^{2\alpha t}\cos ^{2}(\beta t)dt = \dfrac{e^{-2\alpha H}(e^{2\alpha T}-1)}{4\alpha} f(\beta),
\end{split}
\end{equation*}
where
\begin{equation*}
    \begin{split}
  f(\beta) &=  \cos ^{2}(\beta H)+\dfrac{(\alpha \cos (\beta H)-\beta \sin (\beta H))^{2}}{\alpha ^{2}+\beta ^{2}} \\
            &= 1 + \dfrac{\alpha}{ \sqrt{ \alpha^{2} + \beta^{2} }  } \dfrac{\alpha \cos(2\beta H) - \beta \sin(2\beta H)}{ \sqrt{ \alpha^{2} + \beta^{2} } },
  \end{split}
\end{equation*}
hence
\begin{equation*}
v_{1}(\bar{x}_{0}) \leq -\dfrac{\lambda_{min}(W)e^{-2 \alpha H}}{4 \alpha}f(\beta).
\end{equation*}
Notice that, with $\alpha \leq M_{1}$, and $\forall$ $\beta \geq 0$, we have
\begin{equation*}
    \begin{split}
        f(\beta) & \geq \cos ^{2}(\beta H), \\
        f(\beta) & \geq  1-\dfrac{\alpha}{\sqrt{\alpha^{2}+\beta^{2}}} \geq 1 - \dfrac{M_{1}}{\sqrt{M_{1}^{2}+\beta^{2}}}.
    \end{split}
\end{equation*}
Solving the following optimization problem, we can find a positive lower bound (independent of $\beta$) for the function $f(\beta)$:
\begin{equation*}
\underset{\beta \in \lbrack 0,\infty )}{\text{min}}\text{max}%
\left\{ \cos ^{2}(\beta H),1-\dfrac{M_{1}}{\sqrt{M_{1}^{2}+\beta ^{2}}}\right\} .
\end{equation*}
The minimum value of the previous problem is achieved at the point $\beta_{0} \in \left(0,\frac{\pi}{2H}\right)$, such that
\begin{equation*}
\cos ^{2}(\beta _{0} H )=1-\dfrac{M_{1}}{\sqrt{M_{1}^{2}+\beta _{0}^{2}}}.
\end{equation*}
Equation \eqref{Eq. Expression g(b)} is obtained by setting $b_{1}=\beta_{0}H$, and rewriting the resulting terms. Since \eqref{Eq. Expression g(b)} is an increasing function ($g(0)<0$ and $g\left( \frac{\pi}{2} \right)>0$), the number $b_{1}$ on $(0,\pi/2)$, such that $g(b_{1})=0$ exists, and is unique. Therefore, by the fact that $0 < \alpha \leq M_{1}$, we get the desired result expressed by \eqref{Eq. v1<alpha1 varphi en S}.
\end{proof}


\section{Main result}
\label{sec:Main result}
The results of the previous section allow the presentation of two stability criteria for systems with pointwise and distributed delays, that can be verified in a finite number of mathematical operations. 

\subsection{Finite criterion in terms of the Lyapunov matrix}
The first stability result depends exclusively on the delay Lyapunov matrix. 
\begin{theorem}\label{Theorem Stability criterion by finite number of mathematical operation}
System \eqref{Eq. System with distributed delay} is exponentially stable if and only if the Lyapunov condition and the following hold
\begin{equation}\label{Eq. Condition Kr>0}
\mathcal{K}_{\hat{r}} > 0,
\end{equation}
where
\begin{equation}\label{Eq. hat(r)}
\hat{r} = 1 + \lceil e^{LH}H(M_{1}+L)(\alpha^{*}+\sqrt{\alpha^{*}(\alpha^{*}+1)}) -LH\rceil,
\end{equation} 
with $\alpha^{*}=\frac{\alpha_{2}}{\alpha_{1}}$. Here $\alpha_{1}$ and $\alpha_{2}$ are determinated by Theorem \ref{Theorem Stability criterion v1<-a1 } and Lemma \ref{Lemma |v1|<alpha2 |z|<alpha2}, respectively. 
\end{theorem}
\begin{proof}
The necessity directly follows from Theorem \ref{Theorem Criterio Kr}, since condition \eqref{Eq. Condition Kr>0} holds for every number $r$.\\
By contradiction, sufficiency is demonstrated. We assume that system \eqref{Eq. System with distributed delay} is unstable, but $\mathcal{K}_{\hat{r}}>0$ and the Lyapunov condition hold. Therefore, the characteristic equation of the system has no roots on the imaginary axis, and Theorem \ref{Theorem Lyapunov Condition} guarantees the existence and uniqueness of the delay Lyapunov matrix. Consider the function $\hat{\varphi} \in \mathcal{S}$, and $E_{r} = \hat{\varphi} - \psi_{r}$, then
\begin{equation*}
\begin{split}
v_{1}(\psi_{r}) &= v_{1} (\hat{\varphi} - E_{r}) = z(\hat{\varphi} - E_{r},\hat{\varphi} - E_{r}) \\
&= v_{1}(\hat{\varphi}) - 2z(\hat{\varphi},E_{r}) + v_{1}(E_{r}),
\end{split}
\end{equation*} 
where the bilinear function $z(\cdot,\cdot)$ is given by \eqref{Eq. Bilinear functional z}. By Lemmas \ref{Lemma |v1|<alpha2 |z|<alpha2}, \ref{Lemma |Er|<epsilon_r} and Theorem \ref{Theorem Stability criterion v1<-a1 },
\begin{equation*}
\begin{split}
v_{1}(\psi_{r}) &\leq -\alpha_{1}+2\alpha_{2}\|\hat{\varphi}\|_{H}\|E_{r}\|_{H}+\alpha_{2}\|E_{r}\|^{2}_{H}\\
&\leq -\alpha_{1} (1-2\alpha^{*}\varepsilon_{r}-\alpha^{*}\varepsilon_{r}^{2}).
\end{split}
\end{equation*}
By considering $r=\hat{r}$, we have that
\begin{equation*}
1-2\alpha^{*}\varepsilon_{\hat{r}}-\alpha^{*}\varepsilon_{\hat{r}}^{2} \geq 0.
\end{equation*}
Indeed, for $\hat{r}$,
\begin{equation*}
\varepsilon_{\hat{r}} = H\dfrac{(M_{1}+L)e^{LH}}{ \hat{r}-1 + LH}  \leq \dfrac{1}{\alpha^{*}+\sqrt{\alpha^{*}(\alpha^{*}+1)}}
\end{equation*}
and 
\begin{equation*}
\begin{split}
1-2\alpha^{*}\varepsilon_{\hat{r}}-\alpha^{*}\varepsilon_{\hat{r}}^{2} &\geq 1-2\alpha^{*} \left( \alpha^{*}+\sqrt{\alpha^{*}(\alpha^{*}+1)} \right)^{-1} \\ &-\alpha^{*}\left( \alpha^{*}+\sqrt{\alpha^{*}(\alpha^{*}+1)} \right)^{-2}=0,
\end{split}
\end{equation*}
therefore,
\begin{equation*}
  \hat{r} \geq 1+ H(M_{1}+L)e^{LH} \left( \alpha^{*}+\sqrt{\alpha^{*}(\alpha^{*}+1)} \right) - LH. 
\end{equation*}
Finally, from the previous inequality and \eqref{Eq. v1(psi_r)}, we obtain
\begin{equation*}
v_{1}(\psi_{\hat{r}}) = \gamma^{T}\mathcal{K}_{\hat{r}}\gamma \leq 0,
\end{equation*}
which contradicts that $\mathcal{K}_{\hat{r}}>0$ for an unstable system.
\end{proof}
%

\subsection{Finite criterion in terms of the fundamental matrix and the Lyapunov matrix}
The aim of the second criterion is to reduce the number $r$ for which sufficiency is established. The cost of this improvement is the dependence on the system's fundamental matrix. 
Consider the block-matrix
\begin{equation*}
    \mathcal{P}_{r} = \left( I_{n}, K(\delta_{r}),\dots,K((r-1)\delta_{r}) \right) \in \mathbb{R}^{n \times n r}.
\end{equation*}

\begin{theorem}\label{theorem:finite_criterion 2 Kr-alpha0PPT}
System \eqref{Eq. System with distributed delay} is exponentially stable if and only if the Lyapunov conditions holds and
\begin{equation}\label{Eq. Kr-alpha0PP Matrix}
    \mathcal{K}_{r^{*}} - \alpha_{0}\mathcal{P}_{r^{*}}^{T}\mathcal{P}_{r^{*}}>0,
\end{equation}
where
\begin{equation}\label{Eq. r*}
r^{*} = 1 + \lceil e^{LH}H(M_{1}+L)(\alpha^{*}+\sqrt{\alpha^{*}(\alpha^{*}+1)}) -LH\rceil
\end{equation}
with $\alpha^{*}=\frac{\alpha_{2}}{\alpha_{1}+\alpha_{0}}$. Here $\alpha_{0}$, $\alpha_{1}$ and $\alpha_{2}$ are defined in Theorem \ref{Theorem: v1>=alpha0estrella}, Theorem \ref{Theorem Stability criterion v1<-a1 } and Lemma \ref{Lemma |v1|<alpha2 |z|<alpha2}, respectively. 
\end{theorem}
\begin{proof}
\textit{Necessity:} By using the lower bound \eqref{Eq. lower bound v1>=alpha0*} and equation \eqref{Eq. v1(psi_r)}, we have, for every $\gamma \in \mathbb{R}^{n r}$ such that $\psi_r(0)\neq 0$,
\begin{equation*}
	\begin{split}
		\gamma^{T} \left(\mathcal{K}_{r} -\alpha_0 \mathcal{P}_{r}^{T} \mathcal{P}_{r} \right)\gamma
        &=v_1(\psi_{r})-\alpha_0\|\psi_{r}(0)\|^2 \\
        &> v_{1}(\psi_{r}) - \alpha_{0}^{\star}\|\psi_{r}(0)\|^{2} \geq 0.
	\end{split}
\end{equation*}
By Theorem \ref{Theorem Criterio Kr} for the case $\varphi_{r}(0)=0$, $\gamma \neq 0$, the inequality
\begin{equation*}
	\begin{split}
		\gamma^{T} \left(\mathcal{K}_{r} -\alpha_0 \mathcal{P}_{r}^{T} \mathcal{P}_{r} \right)\gamma \geq 0
	\end{split}
\end{equation*}
remains valid.\\
\textit{Sufficiency:} As in the proof of Theorem \ref{Theorem Stability criterion by finite number of mathematical operation},
\begin{equation*}
\begin{split}
v_{1}(\psi_{r}) &= v_{1}(\hat{\varphi}) - 2z(\hat{\varphi},E_{r}) + v_{1}(E_{r})\\
&\leq -\alpha_{1} + 2\alpha_{2}\|E_{r}\|_{H} + \alpha_{2}\|E_{r}\|^{2}_{H}.\\
\end{split}
\end{equation*}
Since $\| \psi_{r}(0) \| = \| \hat{\varphi}(0) \| = 1$,
\begin{equation*}
	\begin{split}
		\gamma^{T} &\left(\mathcal{K}_{r} -\alpha_0 \mathcal{P}_{r}^{T} \mathcal{P}_{r} \right)\gamma
        =v_1(\psi_{r})-\alpha_0\|\psi_{r}(0)\|^2 \\
        &\leq (\alpha_{0}+\alpha_{1}) \left( -1 + 2 \dfrac{\alpha_{2}}{\alpha_{0}+\alpha_{1}} \|E_{r}\|_{H} + \dfrac{\alpha_{2}}{\alpha_{0}+\alpha_{1}} \|E_{r}\|^{2}_{H}\right) \\
        &\leq (\alpha_{0}+\alpha_{1}) \left( -1 + 2 \alpha^{*} \|E_{r}\|_{H} + \alpha^{*} \|E_{r}\|^{2}_{H}\right) \\
	\end{split}
\end{equation*}
The rest of the proof is similar to the one of Theorem \ref{Theorem Stability criterion by finite number of mathematical operation}.\\
\end{proof}
Theorem \ref{Theorem Stability criterion by finite number of mathematical operation} and Theorem \ref{theorem:finite_criterion 2 Kr-alpha0PPT} provide necessary and sufficient conditions for the stability analysis of systems of the form (1), allowing to obtain the stability regions in the space of parameters, including the delays. In both theorems, the value of $r$ depends on the parameters of the system, including the delay, through the bounds of $v_{1}$, the bound of the dynamics of the fundamental matrix of the system as well as of the bound of the derivative of $\varphi$. Therefore, increasing the value of the delay or the system’s parameters results in larger values of $\hat{r}$ and $r^{*}$. Finally, Theorem \ref{theorem:finite_criterion 2 Kr-alpha0PPT} allows reducing the estimate of $r$. However, the conditions depend not only on the Lyapunov matrix but also on the fundamental matrix of the system, which adds to the computational burden as this matrix must be calculated.


\section{Numerical Examples}
\label{sec: examples}

In this section, the stability criteria of Theorem \ref{Theorem Stability criterion by finite number of mathematical operation} and Theorem \ref{theorem:finite_criterion 2 Kr-alpha0PPT} are illustrated by two examples. The implementation is performed in MATLAB, with $\eta=2$ in \eqref{Eq. value of alpha0*}. The Lyapunov matrix associated with the matrix $W=I_{n}$, is computed using the semianalytic method \cite{aliseyko2019lyapunov}. The positivity of $\mathcal{K}_{\hat{r}}$  and ${K}_{r^{*}}-\alpha_{0}\mathcal{P}_{r^{*}}^{T}\mathcal{P}_{r^{*}}$ is checked with the \textit{chol} function, while the function \textit{fzero} is used to find the solution of \eqref{Eq. Expression g(b)}. The fundamental matrix of the system is constructed step by step on the segment $[0,H]$. 
For each example, the values of $\hat{r}$ and $r^{*}$ are calculated for the system parameters that satisfy the necessary conditions of Theorems \ref{Theorem Stability criterion by finite number of mathematical operation} and \ref{theorem:finite_criterion 2 Kr-alpha0PPT} for $r=3$. The stability/instability boundaries obtained by the D-Partitions method are superimposed on the figures.
The numerical computations were performed in a Lenovo Y530 with Intel Core i7-8750H, 2.7 GHz, 6 cores, and 16 GB RAM processor.

\begin{example}\label{Example 1}
In \cite{juarez2020stability}, the stability of a chain of three vehicles, considering human driver's memory effects as distributed delays, is studied. This system is defined by
\begin{equation}\label{Eq: System Example 1}
\dot{e}(t) = A_{1} e(t-h) + \dfrac{1}{2h}G \int_{-3h}^{-h}  e(t+\theta) d\theta,
\end{equation}
with $h = 0.05$ and
\begin{equation*}
A_{1} = \begin{pmatrix}
0 & -2.5\\
-2.5 & 0
\end{pmatrix}, \quad
G = \begin{pmatrix}
-k_{1} & 0 \\
0  & -k_{2}
\end{pmatrix}.
\end{equation*}
The variables $k_{1}$ and $k_{2}$ are design parameters. The integral term defines the maximum delay $H$ being equal to $3h$. Figures \ref{fig:Example 1 r th1} and \ref{fig:Example 1 r th2} represent the maps of the orders of $\hat{r}$ and $r^{*}$ for which the necessary conditions of Theorems \ref{Theorem Stability criterion by finite number of mathematical operation} and \ref{theorem:finite_criterion 2 Kr-alpha0PPT} become sufficient, respectively.  
\begin{figure}
    \centering
    \includegraphics[width=\columnwidth]{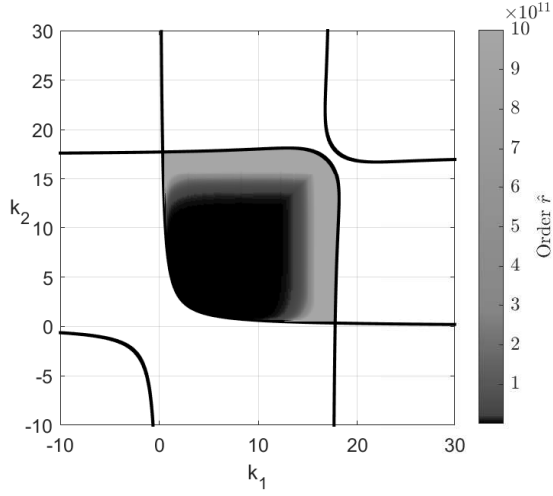}
    \caption{Example \ref{Example 1}: Required orders of $\hat{r}$  in Theorem \ref{Theorem Stability criterion by finite number of mathematical operation} with respect to $(k_{1},k_{2})$.}
    \label{fig:Example 1 r th1}
\end{figure}
\begin{figure}
    \centering
    \includegraphics[width=\columnwidth]{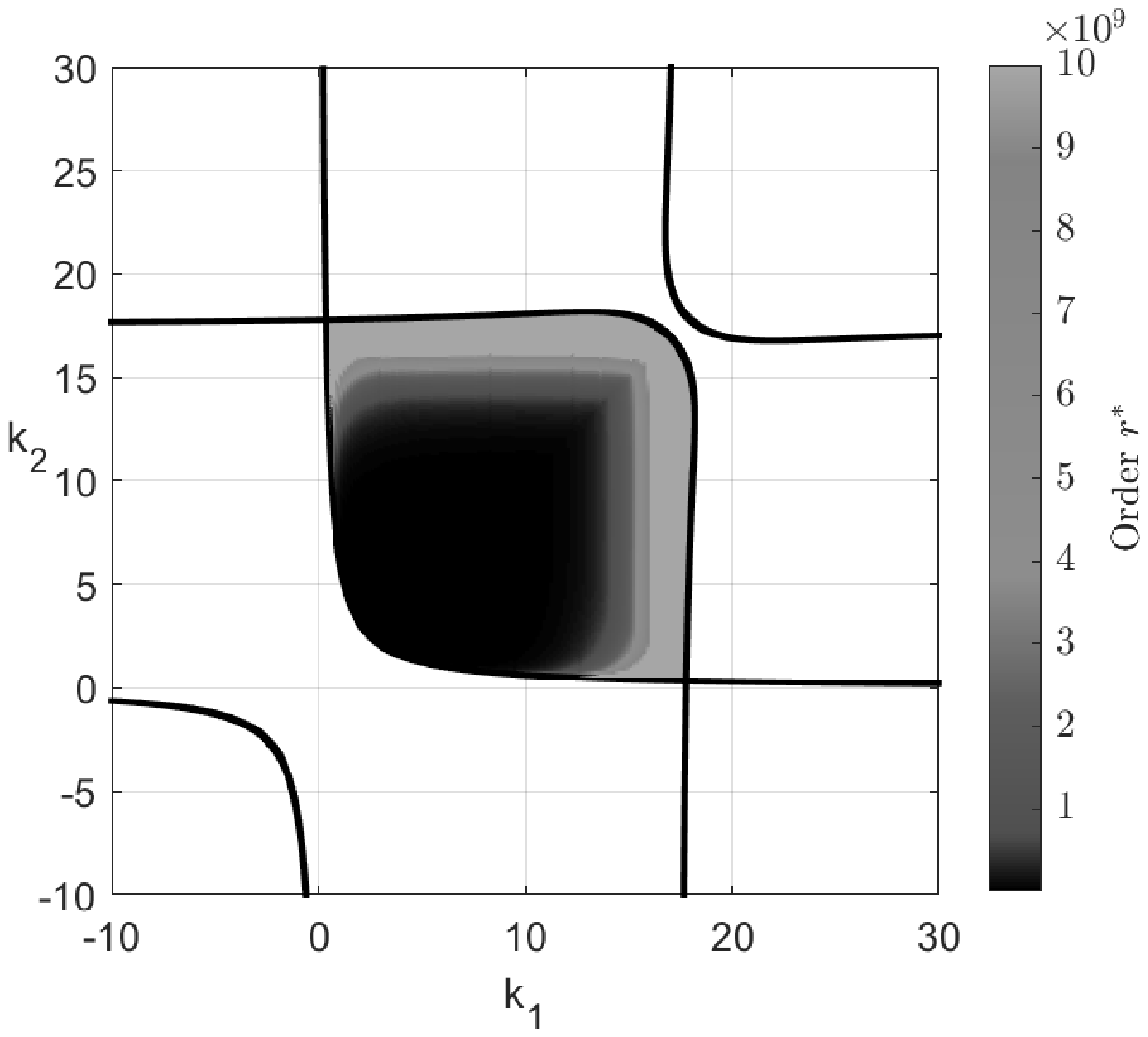}
    \caption{Example \ref{Example 1}: Required orders of $r^{*}$ in Theorem \ref{theorem:finite_criterion 2 Kr-alpha0PPT} with respect to $(k_{1},k_{2})$.}
    \label{fig:Example 1 r th2}
\end{figure}
\end{example}

\begin{example}\label{Example 3}
Consider the following example:
\begin{equation}\label{Eq: System Example 3}
    \dot{x}(t) = \begin{pmatrix}
    -1 & 0.5\\
     0 & p
    \end{pmatrix}x(t-h) + 
    \begin{pmatrix}
     0 & 0\\
    -1 & 0
    \end{pmatrix}  \int_{-h}^{0}x(t+\theta)d\theta
\end{equation}
where $p \in \mathbb{R}$ and $h>0$ are free parameters.
\begin{figure}
    \centering
    \includegraphics[width=\columnwidth]{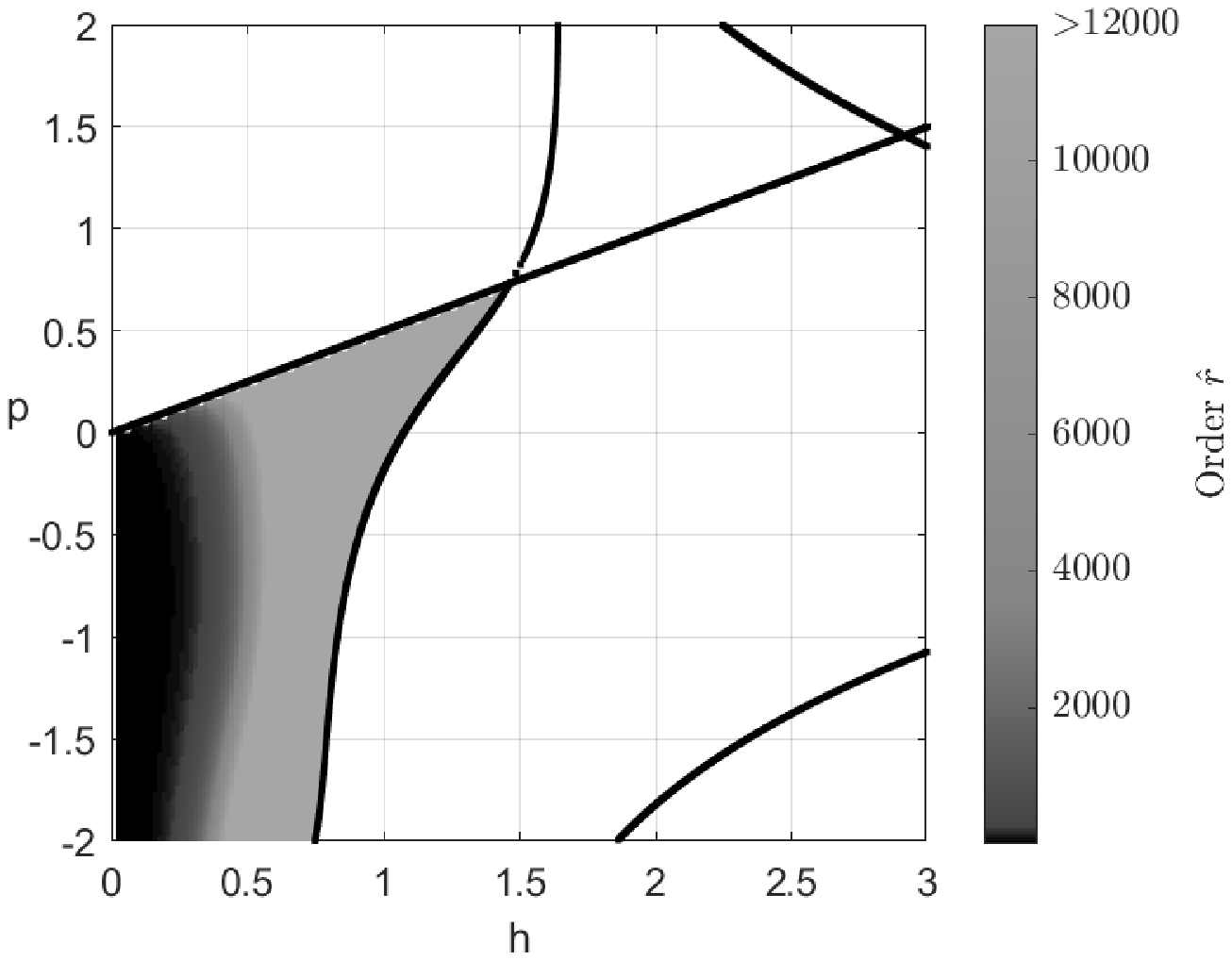}
    \caption{Example \ref{Example 3}: Required orders of $\hat{r}$ in Theorem \ref{Theorem Stability criterion by finite number of mathematical operation} with respect to $(h,p)$.}
    \label{fig:Example 3 r1}
\end{figure}
\begin{figure}
    \centering
    \includegraphics[width=\columnwidth]{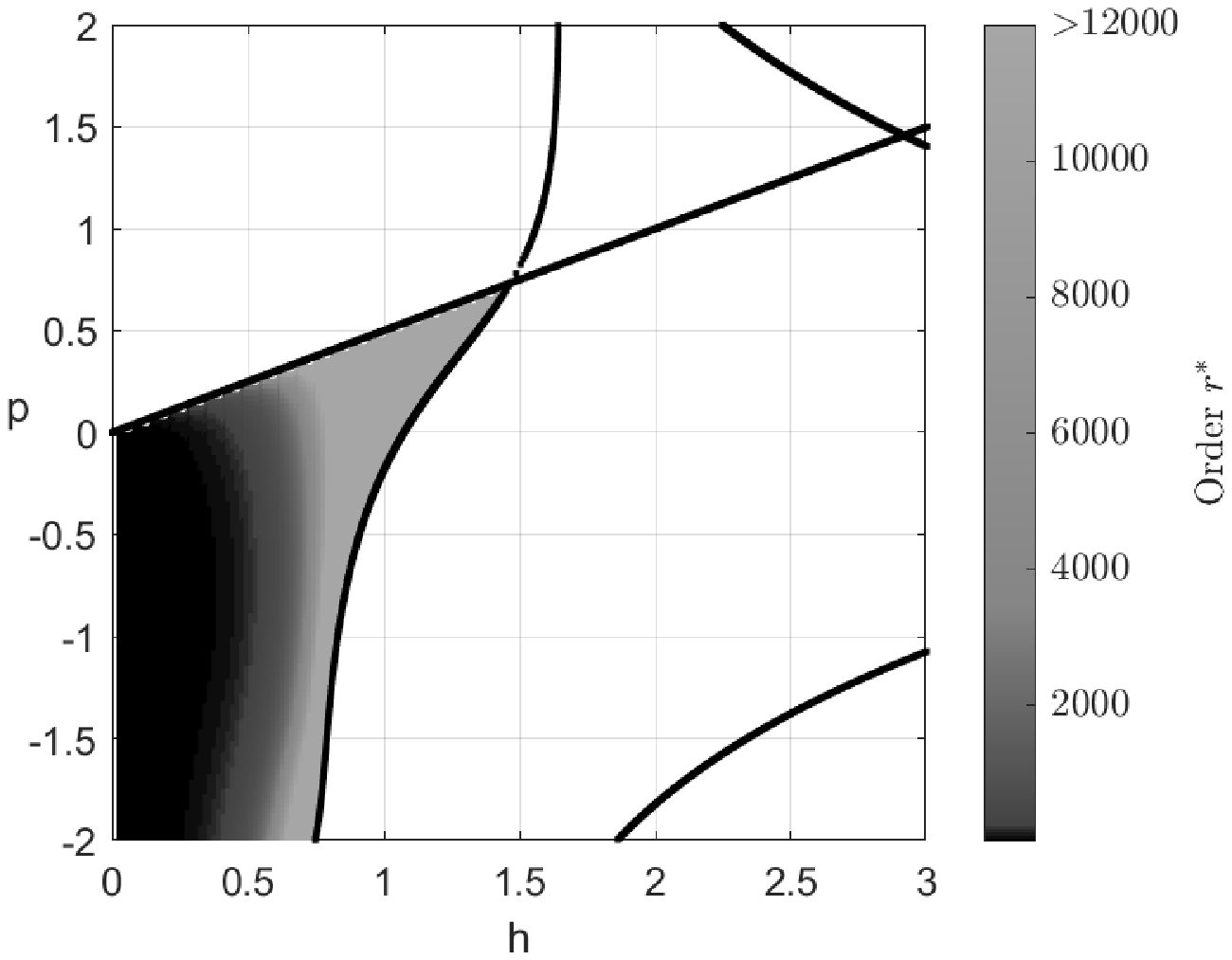}
    \caption{Example \ref{Example 3}: Required orders of $r^{*}$ in Theorem \ref{theorem:finite_criterion 2 Kr-alpha0PPT}  with respect to $(h,p)$.}
    \label{fig:Example 3 r*}
\end{figure}
\begin{table*}
\caption{Experimental results obtained by applying the finite criteria for some parameters $(h,p)$ of Example \ref{Example 3}}
\begin{center}
\begin{tabular}{||c|c|c|c||c|c|c||}
\cline{2-7}
\multicolumn{1}{c|}{} & \multicolumn{3}{||c||}{\textbf{Theorem \ref{Theorem Stability criterion by finite number of mathematical operation}}} & \multicolumn{3}{c||}{\textbf{Theorem \ref{theorem:finite_criterion 2 Kr-alpha0PPT}} }                                                                                     \\ \hline
\multicolumn{1}{||c|}{\textrm{Parameters $(h,p)$}} & \multicolumn{1}{||c|}{$\hat{r}$} & \multicolumn{1}{c|}{Computation time {[}sec{]}} & \multicolumn{1}{c||}{Test result} & \multicolumn{1}{c|}{$r^{*}$} & \multicolumn{1}{c|}{Computation time {[}sec{]}} & \multicolumn{1}{c||}{Test result} \\ \hline
$(0.1,-0.1)$  &	\multicolumn{1}{||c|}{26}  & 0.019  & $\mathcal{K}_{26} > 0$ & 12 & 0.021 & $\mathcal{K}_{12}-\alpha_0\mathcal{P}_{12}^T\mathcal{P}_{12}>0$\\
$(0.25,-0.8)$ &    \multicolumn{1}{||c|}{79}  & 0.026  & $\mathcal{K}_{79} > 0$ & 22 & 0.029 & $\mathcal{K}_{22}-\alpha_0\mathcal{P}_{22}^T\mathcal{P}_{22}>0$\\ 
$(0.3,0.1)$   &	\multicolumn{1}{||c|}{2020}  & 0.338 & $\mathcal{K}_{2020} > 0$ & 463 &  0.068 & $\mathcal{K}_{463}-\alpha_0\mathcal{P}_{463}^T\mathcal{P}_{463}>0$\\
%
%
$(0.2,2)$  &	\multicolumn{1}{||c|}{507}  & 0.058 & $\mathcal{K}_{507} \not\geqslant 0$ & 111 & 0.035 & $\mathcal{K}_{111}-\alpha_0\mathcal{P}_{111}^T\mathcal{P}_{111} \not\geqslant 0$\\
$(0.5,0.5)$ &	\multicolumn{1}{||c|}{9742}  & 4.257 & $\mathcal{K}_{9742} \not\geqslant 0$ & 795 & 0.126 & $\mathcal{K}_{795}-\alpha_0\mathcal{P}_{795}^T\mathcal{P}_{795} \not\geqslant 0$\\
\hline
\end{tabular}
\label{tab:comp_r Example 3}
\end{center}
\end{table*}
The numbers $\hat{r}$ and $r^{*}$ given by \eqref{Eq. hat(r)} and \eqref{Eq. r*} are depicted on Figures \ref{fig:Example 3 r1} and \ref{fig:Example 3 r*}, respectively, for pairs $(h,p)$ in the space of parameters. Table \ref{tab:comp_r Example 3} present a summary of the experiments for selected pairs of parameters $(h,p)$ for Theorem \ref{Theorem Stability criterion by finite number of mathematical operation} (column 2 to 4), and Theorem \ref{theorem:finite_criterion 2 Kr-alpha0PPT}(column 5 to 7). For each result, the estimate of $r$ for which sufficiency holds, the computational time, and the outcome of the test are displayed.
\end{example}


\section{Conclusion} \label{sec:Conclusion}
We have introduced two stability criteria for systems with multiple pointwise and distributed delays. The obtained necessary and sufficient conditions are based on an instability condition inspired by the ideas in \cite{medvedeva2013constructive} and consist of calculating an upper bound of a particular functional. With this result, we can unveil a relation between the stability of system \eqref{Eq. System with distributed delay} and the dimension of \eqref{Eq. Condition Kr>0} and \eqref{Eq. Kr-alpha0PP Matrix}. This fulfills the objective of obtaining a criterion that depends exclusively on the Lyapunov matrix for systems with distributed delays. However, because of the conservatism in the theoretical estimation of $\hat{r}$ the numerical implementation of the stability test demands a high computational effort. The  second criterion reduces the estimated $r^{*}$; however, this result does not depend uniquely on the Lyapunov matrix but also on the fundamental matrix of the system.


\bibliographystyle{IEEEtran}
\bibliography{CNS_Distribuido}

\end{document}